\documentstyle{amsppt}
\magnification=1100
\input xy
\xyoption{all}

\topmatter
\title Unipotent Reduction and the Poincar\'e Problem  \endtitle
\rightheadtext{On the Poincar\'e Problem} \subjclass 14D06, 37F75
\endsubjclass
\author Alexis G. Zamora \endauthor
\address CIMAT,\newline Callej\'on de Xalisco s/n
\newline Valenciana, Guanajuato, Gto.,\newline C.P. 36240, \newline M\'exico \endaddress
\email alexis\@ cimat.mx \endemail
\endtopmatter

\document

\head 1. Introduction \endhead

Let $S$ be a complex projective surface. A foliation on $S$ is an
exact sequence:

$$0@>>>\Cal L @>\Cal F>> TS @>>> \Cal N_{\Cal F}@>>> \Cal
O_{\Gamma}@>>> 0, \tag 1.1 $$ where $\Cal L$ and $\Cal N_{\Cal F}$
are invertible sheaves and $\Gamma$ is a zero-dimensional subscheme
of $S$, called the singular locus of $\Cal F$.

In fact, the entire exact sequence 1.1 is determined by the
inclusion $\Cal F$. Indeed, we can consider the quotient of $\Cal F$
to be:

$$0@>>>\Cal L @>\Cal F>> TS @>>> \Cal C@>>>0.$$ The sheaf $\Cal C$ will not be
locally free in general, but it will be the case for its double dual
sheaf $\Cal N_{\Cal F}:=\Cal C^{**}$. Completing the exact sequence

$$0@>>>\Cal L @>\Cal F>> TS @>>> \Cal N_{\Cal F},$$ we obtain 1.1 ([4], [14]).

The aim of this paper is to give some results concerning the so
called Poincar\'e Problem. Roughly speaking, this problem asks for a
numerical criteria to identify when a morphism $\Cal F:\Cal L@>>>
TS$, defining a foliation, is equal to another of the form $T_f@>>>
TS$, with $f:S@>>> \Bbb P^1$ a holomorphic map and $T_f$ (that must
be isomorphic to our original $\Cal L$) the sheaf of relative vector
fields of $f$.

The original statement of the problem is on $\Bbb P^2$ ([11], [12]),
in this case $f$ must be a rational map $f:\Bbb P^2 \dasharrow \Bbb
P^1$, undefined in a finite number of points (Bezout's Theorem).
Section 2 of this paper is devoted to an explanation of how this
situation can be modified to the case of a holomorphic map $f:S @>>>
\Bbb P^1$, ($S$ will be the blowing-up of $S$ in the indetermination
locus of the original rational map). The classical formulation of
the Poincar\'e Problem is explained there.  This expository section
(which includes, moreover, several results on foliation theory used
below) does not contains any original result and is intended as an
effort to fill a hypothetical gap between the specialists in
foliation theory and those in fibration theory. Once the bridge
between both theories is constructed we work almost completely with
the language of fibration theory.

Section 3 studies the problem of bounding the genus of a fibration:

$$f: S@>>> \Bbb P^1,$$ in terms of numerical information of $K_{\Cal
F}$, the canonical sheaf of the associated foliation. In the
notation introduced in 1.1, $K_{\Cal F}=\Cal L^{-1}(\simeq
T_f^{-1})$. We apply the technique of unipotent reduction of a
fibration to obtain the bounds. In the results proved in that
section some hypothesis on the asymptotical behavior of the
canonical sheaf with respect to the fibers $C$ is made. In this
sense, they are not authentic solutions to the Poincar\'e Problem,
as they assume some knowledge about the fibers $C$. This is an
example of the kind of results we prove in section 3:

\proclaim{Theorem 3.1} Assume $K_{\Cal F}$ is a nef divisor. If, for
some $N>>0$ and some $n\in\Bbb N$

 $$h^0(K_{\bar f}^n)> h^0(K_{\Cal F}^n),$$ then $h^0(K_{\Cal
 F}^n(-C))\ne 0$, and consequently

 $$nK_{\Cal F}^2\ge 2 (g-1).$$
\endproclaim

In this statement $K_{\bar f}^n=\alpha^* K_{\Cal F}^n$, $\alpha$
standing for the unipotent reduction of degree $N$ (see section
3.1).

Our purpose is twofold: on one hand, to show how positive conditions
on $K_{\Cal F}$ (i.e. nefness or ampleness) make the problem easier
to handle; on the other hand, to put the general problem (under the
positive assumptions) within the framework of a well developed
theory.

Finally, section 4 is of a more elementary nature and is devoted to
prove the following result:

\proclaim{Theorem 4.6} Let $\Cal F$ be a foliation defined on $S$
admitting a holomorphic first integral $f:S@>>>\Bbb P^1$ of genus
$g>0$. Assume $\Cal F$ is non-degenerated. If for every singular
point of $\Cal F$, the eigenvalues of the linear part of $\Cal F$
has absolute value greater than $3$ ; then,

$$(g-1)\le 4(-\chi(K_{\Cal F})+\chi (\Cal O_S)).$$ \endproclaim

We prove a number of auxiliary results, each with its own
independent interest. Almost all the arguments used in the proofs
are really simple, but in a crucial step (Lemma 4.4) we use the
unipotent reduction.

This Theorem gives us a complete answer to Poincar\'e Problem when
the eigenvalues associated to the singularities of the reduced
foliation of $\Cal F$ has absolute values greater that 3.

The author wants to express his gratitude  to J. Vitorio Pereira for
several useful comments.

 \head 2.From Poincare's Problem to fibrations
\endhead

Let $\Cal F$ be a foliation on $\Bbb P^2$:

$$0@>>> \Cal O(-m+1) @>\Cal F>> T\Bbb P^2.$$

Assume $\Cal F$ admits a rational first integral. This means that
there exists a rational map:

$$f: \Bbb P^2 \dasharrow \Bbb P^1,$$ such that the set-theoretical fibers
of $f$ coincide with the leaves of $\Cal F$ (in other words: $\Cal
F$ is algebraically integrable).

Rational maps of this kind are given by pencils of curves $\lambda F
+ \mu G=0$ ($\deg F=\deg G=d$).

The classical Poincar\'e Problem consists of finding a function $g$
such that $d\le g(m)$, i.e., to bound the degree of the pencil
defining the first integral in terms of the degree of the foliation.
([11], [12], [15]).

The following basic formula was proved by Poincar\'e. The generic
element of $\lambda F + \mu G=0$ is an irreducible curve. If we
denote by $\sum n_{ij} F_{ij}$ the divisor of reducible components
of the pencil, then:

$$ 2d-2 =m + \sum d_{ij}(n_{ij}-1), \quad d_{ij}:= \deg F_{ij}.$$
Thus, Poincar\'e's Problem is equivalent to bounding the numbers
$n_{ij}$ and $d_{ij}$. In section 4 we deal with the problem of
bounding the multiplicities $n_{ij}$.

A more general (and reasonable) version of the problem is to bound
$d$ in terms of some numerical information of the foliation, and not
only its degree. ([6], [15]).

The aim of this section is to explain how this problem can be
translated into the language of fibered surfaces. In what follows,
we assume that the singularities of $\Cal F$ are non-degenerated.
This implies that the singular points of $\Cal F$ are of two
different classes:

a) Dicritical points: points in the clousure of infinitely many
leaves of $\Cal F$.

b) Saddle points: points in the clousure of exactly two leaves.

Taking into account that these singularities are points where the
map $f$ is not smooth, we obtain the following relations:

a') The dicritical points of $\Cal F$ are the points where the map
$f$ is not defined; equivalently, they are the base locus of the
pencil associated to $f$. The local expression of the pencil
associated to $f$ around these points is of the form: $x^p-\lambda
y^q=0$. The linear part of the local vector field defining $\Cal F$
around a dicritical point is related to the corresponding local form
of $f$: if $f$ is locally $x^p-\lambda y^q=0$, the linear part will
be $diag (p', q')$ with $p=ap'$, $q=aq'$ and $a=m.c.d(p,q).$

b') The saddle points of $\Cal F$ are the singularities of the
non-generic fibers of $f$ (considered with its reduced structures)
located away from the base locus. The local expression  of the
pencil associated to $f$ is of the form: $x^py^q= \lambda $. The
linear part of the local vector field defining $\Cal F$ around a
dicritical point is related to the corresponding local form of $f$:
if $f$ is locally $x^py^q= \lambda $, the linear part will be $diag
(-p', q')$ with $p=ap'$, $q=aq'$ and $a=m.c.d(p,q).$

Now, we can proceed with our translation:

{\bf Step 1} There are two parallel resolution theorems:

1.1. We can solve the indeterminacies of $f$ ([3]): There exists a
minimal chain of blowing-ups:

$$\tilde S=S_n@>>> S_{n-1}@>>> ...@>>>\Bbb P^2,$$ and a commutative
diagram:

$$\matrix \tilde S & & \\ \downarrow^{\pi} & \searrow^{\tilde f} & \\
\Bbb P^2 &\overset{f} \to{\dasharrow}  & \Bbb P^1, \endmatrix$$ such
that $\tilde f$ is a holomorphic map (i.e. defined in all the points
of $\tilde S$).

Therefore, $\tilde f:\tilde S@>>> \Bbb P^1,$ is a fibered surface.
The meaning of {\sl fibration}, as used in this paper, is a morphism
$f: S@>>> Y$ between a nonsingular projective surface $S$ and a non
singular projective curve $Y$, with connected fibers. The general
fiber $\tilde C$ is a non-singular curve, birationally equivalent to
$C$, the clousure of the general fiber of $f$. In particular,
$g(C)=g(\tilde C)$.

1.2. We can solve the dicritical singularities of $\Cal F$. This is
a particular case of a celebrated theorem by Seindeberg ([13]). The
result of this process is a minimal chain of blowing-ups:

$$\tilde S=S_n@>>> S_{n-1}@>>> ...@>>>\Bbb P^2,$$ such that the
foliation $\tilde \Cal F=\pi^* \Cal F$ has only saddle
singularities.

We have used the same notation for the resulting morphism $\tilde
S@>>> \Bbb P^1$ of both processes. They are, in fact, the same:

\proclaim{Proposition 2.1}  Let $\Cal F$ be a non-degenerated
foliation on $\Bbb P^2$ admitting a rational first integral

$$f: \Bbb P^2 \dasharrow \Bbb P^1.$$

The resolution process of Seindenberg coincides with the resolution
of indeterminacies of $f$. In particular, $\tilde f:\tilde S@>>>
\Bbb P^1,$ is a holomorphic first integral for $\tilde \Cal F$ and
the geometric genus of the general leaves of $\Cal F$ and $\tilde
\Cal F$ are the same. \endproclaim

The proof of this fact is, under the non-degenerated hypothesis, an
elementary calculus based on the local equations of the blowing-up.

{\bf Step 2.} Poincar\'e  ([11], see also [9] and [15]) proved the
following inequality: if $\Cal F$ is non-degenerated, of degree $d$
and admits a rational first integral; then, assuming our previous
notation:

$$\frac{m-4}{4}d+1\le g. \tag 2.1$$

In virtue of this relation and Step 1, the problem can be
reformulated as: given a fibration $f: S@>>>\Bbb P^1$ of genus $g$,
find an upper bound for $g$ in terms of numerical invariants of the
foliation associated to $f$.

The last phrase needs some explanation. Associated to $f$ there is
an invertible sheaf $K_f:= K_S\otimes f^*K_{\Bbb P^1}^{-1}$, the so
called relative dualizing sheaf. On the other hand, the foliation
associated to $f$ is given by an inclusion:

$$0@>>> \Theta_f @> \Cal F>>TS @>df>> f^*T\Bbb P^1 @>>> \Cal I
@>>>0. $$

The sheaf $\Cal I$ is supported on the singular points of $f$. The
sheaf $(\Theta_f)^{-1}$ is denoted by $K_{\Cal F}$ and called the
{\sl canonical sheaf of $\Cal F$} ([4], [8]).

We introduce the following notation, used systematically in the
remaining part of this paper. Denote by

$$\Delta =\sum_{ij} n_{ij} F_{ij},$$ the divisor defined as the sum of all
the non-reduced fibers of $f$. Hence, if we fix the index $i$,
$\sum_j n_{ij}F_{ij}$ will denote a single non-reduced fiber of $f$.

Moreover, we write $\Delta=\Delta_{red}+\Delta_0$, where
$\Delta_{red}=\sum_{ij} F_{ij}$ is the reduced divisor associated to
$\Delta$. Note that $\Delta_0=\sum_{ij} (n_{ij}-1) F_{ij}$. The
relation between $K_f$ and $K_{\Cal F}$ is:

$$K_{\Cal F}=K_f(-\Delta_0).$$

In [14] the following was proved:

\proclaim{Theorem 2.2} Let $f:S@>>> Y$ be a fibration of genus $g\ge
2$. The following are equivalent:

a) $K_{\Cal F}$ is big and nef,

b) $f$ is relatively minimal (no exceptional curves of the first
kind contained in the fibers) and non-isotrivial (two general fibers
are not isomorphic).

Furthermore, if a) or b) holds, then the only irreducible curves
$D\subset S$ such that $K_{\Cal F}\cdot D=0$ are the $(-2)$-vertical
curves (rational curves with self-intersection $-2$ contained in
some fiber of $f$).
\endproclaim

The proof of this theorem uses the unipotent reduction of $f$. The
next section is devoted to the exploitation of the same technique in
connection with Poincar\'e's problem.

\head 3. Unipotent reduction and the Poincar\'e Problem \endhead

\subhead 3.1 Unipotent reduction of a fibration\endsubhead In this
section we apply a useful technique in fibration theory, namely the
unipotent reduction, to the study of the Poincar\'e Problem. We
first recall the basic principles of unipotent reduction.

Let $f:S @>>> Y$ be a fibration, then there exists a commutative
diagram:

$$ \xymatrix  {\bar S \ar[r]^{\alpha} \ar[d]_{\bar f} & S \ar[d]^f \\
X \ar[r]_{N:1}^{\pi}& Y},  $$ such that $\bar S$ is a non-singular
projective surface and $\bar f$ is a reduced fibration, i.e. all its
fibres are reduced. The ramification values of $\pi$ contain the
critical values of $f$ as a subset.

In this construction $N$ is a common multiple of the multiplicities
of the non reduced components of the fibers of $f$, and $\pi$ is a
totally ramified cyclic covering. We assume that $N$ is sufficiently
large and $\pi$ has $N$ ramification values (see [2], sections III.9
and III.10, for the details concerning this construction). The genus
of $X$ can be computed easily, for sufficiently large $N$, using
Riemann-Hurwitz and the fact that $\pi$ is totally ramified on $N$
points:

$$g_X-1=\frac{N(N-1)}{2}+N(g_Y-1).\tag 3.1$$ In particular, if $Y=\Bbb P^1$,
$g_X-1=\frac{N(N-3)}{2}.$

The other fundamental fact that we shall use sys\-tem\-ati\-cally is
due to Serrano:

$$K_{\bar\Cal F}:=\alpha^* K_{\Cal F}$$ is the relative dualizing
sheaf of $\bar f$. Thus $K_{\bar\Cal F}=K_{\bar f}= K_{\bar
S}\otimes {\bar f}^* {K_X}^{-1}$ ([14], claim in the proof of Prop.
2.1), we shall call this result  Serrano's Lemma. Thus, $K_{\bar
\Cal F}.\bar C=2g-2$ where $\bar C$ denotes a fiber of $\bar f$.

\subhead 3.2 Positivity properties of $K_{\bar\Cal F}$ and explicit
bounds for $g$ \endsubhead

Let $f:S@>>>\Bbb P^1$ be a fibration with fibers supported on normal
crossing divisors (i.e., the associated foliation is
non-degenerated). As explained in section 2, it is an important
problem in foliations theory to find a upper bound for the genus $g$
of the fibres of $f$ in terms of the invariants of $S$ and $K_{\Cal
F}$. The aim of this section is to show how positivity conditions on
$K_{\Cal F}$ allows us to solve this problem.  In the remains of
this paper  we use the notation introduced in section 3.1.

\proclaim{Theorem 3.1} Assume $K_{\Cal F}$ is a nef divisor. If, for
some $N>>0$ and some $n\in\Bbb N$,

 $$h^0(K_{\bar f}^n)> h^0(K_{\Cal F}^n),$$ then $h^0(K_{\Cal
 F}^n(-C))\ne 0$, and consequently

 $$nK_{\Cal F}^2\ge 2 (g-1).$$
\endproclaim

\demo{Proof} The direct image sheaf $\alpha_*\Cal O_{\bar S}$ can be
computed explicitly ([5]):

$$\alpha_*\Cal O_{\bar S}=\oplus_{i=0}^{N-1} \Cal L^{-(i)}.$$ As $f$ takes values in $\Bbb P^1$
the first terms in this decomposition are: $\Cal L^{(0)}=\Cal O_S$,
$\Cal L^{-(1)}=\Cal O_S(-C)$. In general, $\Cal L^{-(i)}$ is an
invertible sheaf contained in $\Cal L^{-(i-1)}$:

$$0@>>> \Cal L^{-(i)}@>>>\Cal L^{-(i-1)}. \tag 3.2 $$

Until this point, all the considerations have been general. Now, in
order to simplify the notation, we develop the complete argument in
the case $n=1$. The general case follows after an obvious
modification.

Using projection formula and Serrano's Lemma we obtain:

$$\align \alpha_*K_{\bar f}&\simeq \alpha_* (\alpha^*K_{\Cal F})=K_{\Cal
F}\otimes \alpha_*\Cal O_{\bar S} \\ &=K_{\Cal
F}\otimes (\Cal O_S\oplus \Cal O_S(-C)\oplus ...)\\
&= K_{\Cal F}\oplus K_{\Cal F}(-C) \oplus ... \endalign $$

If $h^0(K_{\bar f})> h^0(K_{\Cal F}),$ then, by inclusion 3.2,
$h^0(K_{\Cal F}(-C))>0$. The Theorem follows from the fact that
$K_{\Cal F}$ is nef:

$$K_{\Cal F}\cdot(K_{\Cal F}(-C))=K_{\Cal F}^2 -2(g-1)\ge 0.$$

\qed \enddemo

Using Theorem 3.1 we can deduce another condition ensuring the
existence of a bound for $g$:

\proclaim{Theorem 3.2} Assume $K_{\Cal F}$ is an ample line bundle.
Let $n$ be a natural number such that $K_{\Cal F}^n$ is very ample.
If, for $N>>0$, $h^0(K_{\bar f}^n(-C))>0$, then

$$n^2K_{\Cal F}^2\ge 2 (g-1).$$ \endproclaim

\demo{Proof} Once again, we concentrate on the case $n=1$. We can
assume that $h^0(K_{\bar f})= h^0(K_{\Cal F})$, otherwise, we obtain
the conclusion from Theorem 3.1 (obviously, $K_{\Cal F}$ very ample
is more than enough to guarantee it is nef).

Fix the embedding $S@>>> \Bbb P^k$ given by the linear system
$K_{\Cal F}$. If we apply unipotent reduction, we obtain a
commutative diagram:

$$ \xymatrix  {\bar S \ar[r]^{\alpha} \ar[d]_{\vert K_{\bar f}\vert} & S \ar[d]^{\vert  K_{\Cal F}\vert,}\\
S_0 \ar[r]& S} $$(we abuse the notation slightly by using the same
symbol for a linear system and the rational map it defines).

The bottom arrow is simply the projection on the subspace
$\alpha^*\vert K_{\Cal F}\vert$. Under our hypothesis this last
morphism must be the identity and $\alpha=\vert K_{\bar f}\vert$.

Let $C_1$ be a fiber such that $h^0(K_{\bar f}(- C_1))\ne 0$, call
$C_2$, ..., $C_N$, the set of fibers of $\bar f$ such that $\alpha
(C_i)=\alpha(C_1)$. Then, since $\alpha=\vert K_{\bar f}\vert$, we
must have $h^0(K_{\bar \Cal F}(-\bar C_1 -\bar C_2-...-\bar C_N))\ne
0$. Let $\bar B\in \vert K_{\bar f}\vert$ be an  irreducible and
nonsingular curve (the existence of such a curve is guaranteed by
Bertini's Theorem). Restricting a non zero section of $ K_{\bar
f}(-\bar C_1 -\bar C_2-...-\bar C_N)$ to $\bar B$, we obtain a
nonzero section of $K_{\bar f}(-\bar C_1 -\bar C_2-...-\bar
C_N)\vert_{\bar B}$. Thus this sheaf is of positive degree and:

$$\bar B\cdot (\sum_{i=1}^N C_i)=N(2g-2)\le K_{\bar f}\cdot
\bar B=K_{\bar f}^2 \le N K_{\Cal F}^2.$$ \qed

\enddemo

We don't know of any argument ensuring the existence of global
sections of $h^0(K_{\bar f}^n(-C))$. However, we can, at least,
prove the following:

\proclaim{Lemma 3.3} Assume $K_{\Cal F}^n$ is big and nef. Then, for
all $N>>0$ and any fiber $C$ of $\bar f$, $K_{\bar f}^n(-C)$ is big
and nef. \endproclaim

\demo{Proof} This lemma is valid even if the original fibration
$f:S@>>> Y$ is defined over a curve $Y$ of genus greater than zero.
As previously, we make the proof for $n=1$.

Obviously $K_{\bar f}(-\bar C)$ is big for $N$ large enough. Let $D$
be any irreducible curve on $\bar S$. We must to prove $K_{\bar
f}(-C)\cdot D\ge 0$. We can assume $D$ is not contained in a fiber
of $\bar f$, because for any irreducible component $F$ of a fiber,
$F\cdot C=0$ and $K_{\bar f}$ is nef. Let $\alpha (D)=D_0$. Note
that $D_0$ is an irreducible curve in $S$, as all the exceptional
curves of $\alpha$ are contained in fibers of $\bar f$. We have a
commutative diagram:

$$ \xymatrix  {D \ar[r]^{\alpha\vert_D} \ar[d]_{\bar f\vert_D} & D_0 \ar[d]^f\vert_{D_0} \\
X \ar[r]^{\pi}& Y}.  $$

The degree of $\bar f\vert_D$ is $\bar C. F$ and the degree of
$f\vert_{D_0}$ is $C. D_0$, whereas the degree of $\alpha\vert_D$ is
a divisor $n$ of $N$. Thus, $C. D_0=\frac{N}{n}\bar C. D$. Moreover,
$\bar C. D\ge C. D_0$, so we conclude $n=N$ and $\bar C. D= C\cdot
D_0$. From this it follows that $\alpha^*D_0=D+E$ where all the
irreducible components of $E$ are exceptional curves of $\alpha$.

Observe that any irreducible exceptional curve $A$ of $\alpha$ is a
$(-2)$ curve, so $K_{\bar f}. A=0$ as a consequence of Serrano's
Lemma and adjunction formula. It follows that $K_{\bar f}. E=0$.
Applying the projection formula for the intersection of divisors we
get:

$$\align K_{\bar f}(-\bar C).D =& D. K_{\bar f}-D. \bar C \\ & =(D+E). K_{\bar f} -E. K_{\bar f} -D. \bar
C \\ & = \alpha^*D_0 \cdot \alpha^*K_{\Cal F}- D\cdot \bar C
\\ & = NK_{\Cal F}. D_0-C. D_0.\endalign$$

The only curves in $S$ having zero intersection with $K_{\Cal F}$
are $(-2)$ curves contained in the fibers of $f$ (Theorem 2.2).
Therefore, $K_{\Cal F}. D_0>0$, and for $N>>0$:

$$ K_{\bar f}(-\bar C).D>0.$$

\qed
\enddemo

A natural way to study the space of sections $H^0(K_{\bar f}^n(-C))$
is considering its direct image under $\bar f$. We have another
partial result:

\proclaim{Proposition 3.4} Assume $K_{\Cal F}$ is big, nef and
effective. If for some $N>>0$ and $n\ge 2$, the direct image sheaf
$\bar f_*K_{\bar f}^n$ splits in direct sum of invertible sheaves,
then there exists a fiber $C$ such that: $H^0(K_{\bar f}^n(-C))\ne
0$.
\endproclaim

\demo{Proof} Once again, we consider the case of the lowest possible
power, i.e., $n=2$. The proposition is a consequence of the
positivity theorem of Arakelov, Parshin and Fujita ([1], [2]). In
the first place, $K_{\bar f}^2$ will be big and nef. Our assumption
is:

$$\bar f_*K_{\bar f}^2=\oplus_{i=1}^{3g-3}L_i,$$ with $L_i$ invertible.

Moreover, $\deg L_i\ge 0$, for all $i$. Moreover, $L_i\ne \Cal O_X$
for all $i$, because $h^1(-K_{\bar f}^n)=0$ ([1]). Thus, for some
$i$ (say $i=1$), we must have $h^0(L_1)>0$ and $\deg L_1>0$.

Write $L_1=\Cal O_X(p_1+...+p_r)$, then $h^0(L_1(-p_1))\ne 0$. As
$L_i$ is a sumand of $\bar f_*K_{\bar f}^2$ and using the projection
formula, we conclude that $h^0(\bar f_*K_{\bar f}^2(-C_{p_1}))\ne
0$, where $C_{p_1}$ denotes the fiber over $p_1$ . \qed \enddemo

A last application of this circle of ideas is:

\proclaim{Proposition 3.5} If, for some $N>>0$ and a natural number
$n$, the direct image sheaf  $\bar f_*K_{\bar f}^n$ is generated by
global sections, then, either:

$$(2n-1)(g-1)\le h^0(K_{\Cal F}^n) \quad \text{or } 2(g-1) \le nK_{\Cal F}^2.$$

\endproclaim

\demo{Proof} We explain the case $n=1$. The hypothesis means that
there exists a surjective morphism of locally free sheaves:

$$ H^0(\bar f_*K_{\bar f})\otimes \Cal O_X= H^0(K_{\bar f})\otimes \Cal
O_X @>>> \bar f_*K_{\bar f} @>>>0.$$

The rank of $\bar f_*K_{\bar f}$ is $g$ and the rank of $H^0(\bar
f_*K_{\bar f})\otimes \Cal O_X$ is $h^0(K_{\bar f})$. It follows at
once that $g\le h^0(K_{\bar f})$. By Theorem 3.1 we can assume that
$h^0(K_{\bar f})=h^0(K_{\Cal F})$, and we obtain the result.

\qed
\enddemo

\head 4. Bounding the multiplicities of non-generic fibres \endhead

Let us start with a fibration  $f:S @>>> \Bbb P^1$. Associated to
$f$ we have the relative duality sheaf

$$K_f:=K_S\otimes f^*K_{\Bbb P^1}^{-1}=K_S(2C), \tag 4.1$$
since $K_{\Bbb P^1}=\Cal O_{\Bbb P^1}(-2)$. In the same way we have
the canonical sheaf associated to the foliation $\Cal F$ defined by
the fibres of $f$. This sheaf, denoted by $K_{\Cal F}$, is related
to $K_f$ by means of the formula

$$K_f=K_{\Cal F}(\Delta_0)\tag 4.2 $$ (see the notation introduced in section
2).

Our first goal is the following:

\proclaim{Lemma 4.1} Let $f: S@>>>\Bbb P^1$ be a fibration; then:

$$\Delta_0^2= (K_{\Cal F}-K_S)^2.$$
\endproclaim

\demo{Proof} From relations 4.1 and 4.2 we obtain:

$$  K_{\Cal F}-K_S = (K_S+2C -\Delta_0)-K_S=2C-\Delta_0.$$

Since all the irreducible components of $\Delta_0$ are contained in
fibers of $f$, we have $C\cdot \Delta_0=0$. Hence:

$$  (K_{\Cal F}-K_S)^2 = (2C-\Delta_0)^2= \Delta_0^2.$$ \qed
\enddemo

In this way, we have determined the value of $-\Delta_0$ in terms of
numerical information of the foliation $\Cal F$. For what follows,
denote by $s$ the number of non-reduced fibers of $f$. Note that

$$ \Delta_0^2= (\Delta
-\Delta_{red})^2 = (sC-\Delta_{red})^2=\Delta_{red}^2.$$

Now, fix a critical value of $f$, $t_i\in \Bbb P^1$. Let $C_i$ be
the corresponding fibre, writing $C_i=\sum_j n_{ij}F_{ij}$, we have:

$$ (\sum_j F_{ij})^2= \sum_j F_{ij}^2+2 \sum_{j<k} F_{ij}\cdot F_{ik}.$$

We denote by $c_2(\Cal F)$ the number of singular points of $\Cal
F$. Taking into account that the singularities of $\Cal F$ are
precisely the intersections of the irreducible components of the
non-generic fibers of $f$ we obtain:

$$\Delta_0^2=\Delta_{red}^2= \sum_{ij} F_{ij}^2 +2 c_2(\Cal F).$$

Using this equality and Lemma 4.1 we can formulate our main
conclusion to this point:

\proclaim{Proposition 4.2} With the same notation as before we have:

$$-\sum_{ij} F_{ij}^2 =(K_{\Cal F}-K_S)^2+2c_2(\Cal F). $$
\endproclaim

Our next observation is that, if we fix a non-reduced fiber, say
$\sum_jn_{ij}F_{ij}$, and we know a bound for a single multiplicity
$n_{ij_0}$, then we can obtain a bound for the remainder
multiplicities $n_{ij}$. In fact:

$$\align 0= F_{ij_0}\cdot C &= F_{ij_0}\cdot (\sum_j n_{ij}F_{ij})\\ &=
n_{ij_0}F_{ij_0}^2+ \sum_{j\ne j_0} n_{ij}F_{ij_0}\cdot F_{i,j}.
\endalign $$

Thus, the numbers $F_{ij_0}^2$ are computed as:

$$F_{ij_0}^2= -\frac{\sum_{j\ne j_0} n_{i,j}F_{ij_0}\cdot F_{i,j}}{n_{ij_0}}.\tag 4.3$$

By Proposition 4.2, the number $-F_{ij_0}^2$ is bounded by $\Cal F$.
Thus, a bound for $n_{ij_0}$ gives a bound for the multiplicities of
the components intersecting positively to $F_{ij_0}$. Using the fact
that all the fibers of $f$ are 1-connected and iterating this
process we obtain a bound for all the $n_{ij}$ ($i$ fixed).

This observation motivates the principal result of this section:

\proclaim{Theorem 4.3} Let $f:S@>>> \Bbb P^1$ be a fibration of
genus $g\ge 2$, such that all its singular fibers have support on
normal crossing divisors. Let $n_{ij}$ be the multiplicities of the
components of the singular fibers of $f$. Then, either:

. there are at least $s-2$ non-reduced fibers of $f$, such that some
of the associated multiplicities $n_{ij}$ are less than or equal to
$3$, or

. $(g-1)\le 4(-\chi(K_{\Cal F})+\chi (\Cal O_S))$
\endproclaim

Before starting the proof we need some preliminaries. We can apply
the unipotent reduction to $f$ (see section 3). Denote by $\bar
F_{ij}$ the pre-image of the curves $F_{ij}$ under $\alpha$, the
morphism defining the unipotent reduction. With this notation in
mind we have the following lemma:

\proclaim{Lemma 4.4}

$$(s-2)(g-1) +\chi(K_{\Cal F})-\chi(\Cal O_S)= \sum_{ij} \frac{\bar
F_{ij}\cdot K_{\bar f}}{2n_{ij}}.$$ \endproclaim

\demo{Proof} The proof makes use of the unipotent reduction of $f$,
explained in the previous section. Let $\bar f: \bar S@>>> X$ be the
unipotent reduction of $f$ associated to $N$, a common multiple of
$n_{ij}$ for all $i$ and $j$. The number $K_{\bar f}\cdot K_{\bar
S}$ can be computed in two different ways:

$$\align K_{\bar f}\cdot K_{\bar S}=& K_{\bar f}\cdot(K_{\bar f}\otimes\bar
f^*K_X)\\ &= NK_{\Cal F}^2+ 4(g-1)(g_X-1), \endalign$$ and, using
that $\alpha$ is a ramified covering away from some exceptional
curves that are contracted to points, we can use the Riemann-Hurwitz
formula ([2]):

$$\align K_{\bar f}\cdot K_{\bar S} =&\alpha^*K_{\Cal
F}\cdot(\alpha^* K_S + \sum_{ij} (N/n_{ij}-1)\bar F_{ij}+
(N-1)\sum_{k=1}^{N-s}C_k) \\ &= NK_{\Cal F}\cdot K_S+\sum_{ij}
(N/n_{ij}-1)\bar F_{ij}\cdot K_{\bar f}+ 2(N-s)(N-1)(g-1).
\endalign$$

If we consider these equalities as polynomials in $N$, we obtain the
result from the degree 1 term. \qed \enddemo

Now, we can reach the proof of Theorem 4.3:

\demo{Proof} Denote by $l$ the number of non-reduced fibers with
some associated multiplicity less than or equal to $3$. After
reordering the non-reduced fibers, we can assume that these $l$
fibers are exactly $C_1$,...,$C_l$. Assume $l\le s-3$. Using the
previous lemma write:

$$\align (s-2-l)(g-1)+ l(g-1)-\sum_{i=1}^l\frac{\bar F_{ij}\cdot K_{\bar
f}}{2n_{ij}}& =\sum_{i=l+1}^s \frac{\bar F_{ij}\cdot K_{\bar
f}}{2n_{ij}} -(\chi(K_{\Cal F}) -\chi(\Cal O_S))\\ &\le
\frac{(s-l)(g-1)}{4} -(\chi(K_{\Cal F}) -\chi(\Cal O_S)).\endalign$$

The inequality uses the fact that, for any $i$, $(\sum_j \bar
F_{ij})\cdot K_{\bar f}=2(g-1)$, as the first divisor is a fiber of
$\bar f$. For the same reason:

$$l(g-1)-\sum_{i=1}^l\frac{\bar F_{ij}\cdot K_{\bar
f}}{2n_{ij}}\ge 0.$$

Thus, we only need to prove that $$4(s-l-2)-(s-l)\ge 1.$$ This
follows from $l\le s-3$. \qed \enddemo

\proclaim{Corollary 4.5} Using the same notation as in Theorem 4.3
we have, either:

. for at least $s-2$ non-reduced fibers of $f$, its multiplicities
are bounded by $\Cal F$, or

. $(g-1)\le 4(-\chi(K_{\Cal F})+\chi (\Cal O_S))$. \endproclaim

\demo{Proof} Simply combine Theorem 4.3 and the observation
preceding it. \qed \enddemo

Theorem 4.3 allows us to obtain an answer to Poincar\'e's Problem
under some hypotheses on the singular points of $\Cal F$ (compare
with [10], Chapter 7, Theorem 15). We use the terminology explained
in section 2.

\proclaim{Theorem 4.6} Let $\Cal F$ be a foliation defined on $S$
admitting a holomorphic first integral $f:S@>>>\Bbb P^1$ of genus
$g>0$. Assume $\Cal F$ is non-degenerated. If for every singular
point of $\Cal F$, the eigenvalues of the linear part of $\Cal F$
have absolute value greater than $3$ ; then,

$$(g-1)\le 4(-\chi(K_{\Cal F})+\chi (\Cal O_S)).$$ \endproclaim

\demo{Proof} The condition $p,q>3$ for all the singular points
implies that all the multiplicities of non-reduced fibers of $f$
satisfy $n_{ij}>3$ (see item b') in section 2). Thus, the theorem
follows at once from theorem 4.3. \qed \enddemo

\enddocument